\documentclass[11pt,twoside]{article}

\usepackage[ansinew]{inputenc} 
\usepackage{amsmath,amsfonts,amssymb,amsthm}
\usepackage{pstricks,pst-node,pst-coil,pst-plot,pstricks-add}
\usepackage{graphics,geometry,epsfig}
\usepackage{bbm}

\usepackage{mathrsfs} 
\newcommand{\LL}{\mathscr L}
\newcommand{\dom}{D}
\newcommand{\eps}{\varepsilon}
\newcommand{\weak}{\rightharpoonup}
\newcommand{\uA}{\underline A}

\newcommand{\Rea}{\operatorname{Re}}

\newcommand{\N}{\mathbb{N}}
\newcommand{\R}{\mathbb{R}}

\newtheorem{theorem}{Theorem}[section]

\theoremstyle{remark}
\newtheorem{remark}[theorem]{Remark}

\title{Well-posedness of Non-autonomous Linear Evolution Equations in Uniformly Convex Spaces}

\author{ Jochen Schmid and Marcel Griesemer\\  
\small Fachbereich Mathematik, Universit\"at Stuttgart, D-70569 Stuttgart, Germany\\
\small firstname.lastname@mathematik.uni-stuttgart.de}  

\date{}

\begin{document}
\maketitle                 

\begin{abstract}
This paper addresses the problem of well-posedness of non-autonomous
linear evolution equations $\dot x = A(t)x$ in uniformly convex Banach
spaces. We assume that $A(t):D\subset X\to X$ for each $t$ is the
generator of a quasi-contractive, strongly continuous group, where the domain $D$ and the growth
exponent are independent of $t$. Well-posedness holds
provided that $t\mapsto A(t)y$ is Lipschitz for all $y\in D$. H\"older
continuity of degree $\alpha<1$ is not sufficient and the assumption of uniform convexity cannot be dropped.
\end{abstract}

\section{Introduction}
In the literature the existence of the propagator (evolution system) for the non-autono\-mous Schr\"odinger equation
is often discussed within the more general context of abstract non-autonomous linear evolution equations
\begin{equation}\label{IVP}
      \dot x = A(t)x, \qquad x(s)=y
\end{equation}
in some Banach space $X$, where $A(t):D(A(t))\subset X\to X$ for each $t\in
[0,T]$ is the generator of a strongly continuous semigroup, $0\leq s\leq t$ and $y\in D(A(s))$. On the level of proofs this
approach involves serious technical difficulties that are associated
with the lack of structure of general Banach spaces and the
non-reversibility of the dynamics given by a semigroup. The prize for 
the solution of these problems is paid in terms of regularity assumptions
on $t\mapsto A(t)$ \cite{Kato70,Kato73,Dorroh75,Tanabe79,Pazy,Goldstein}. 

In the present paper,  which is motivated by the Schr\"odinger equation, 
the evolution problem \eqref{IVP} is discussed in a more restrictive setting, which does
not have the drawbacks mentioned above. In this setting  
$X$ is a uniformly convex Banach space and $A(t)$ for each $t\in
[0,T]$ is the generator of a strongly continuous \emph{group} rather than a semigroup. We
assume, moreover, that this group is quasi-contractive with a growth exponent that is independent of $t$
and that the domain $D=D(A(t))$ is independent of $t$ as well.
Our main result, in the simplest form, establishes the existence of a
unique evolution system $U(t,s)$ provided 
\begin{equation}\label{curve} 
    t\mapsto A(t)y
\end{equation}
is Lipschitz for all $y\in D$.  It follows that $t\mapsto x(t)=U(t,s)y$ is the unique continuously differentiable
solution of \eqref{IVP} and that it depends continuously on the
initial data $s$ and $y$ (well-posedness). We give examples showing
that H\"older continuity of the map \eqref{curve} is not sufficient
and that Lipschitz continuity is not sufficient anymore if the assumption of 
uniform convexity is dropped. This means in 
particular that well-posedness of the non-autonomous Schr\"odinger equation, that is,
Equation~\eqref{IVP} with $X$ a Hilbert space and
$A(t)^{*}=-A(t)$, requires less regularity than well-posedness of
\eqref{IVP} in the general Banach space setting. 

The well-posedness of \eqref{IVP} in uniformly convex Banach spaces
was previously studied by Kato~\cite{Kato53, Kato70}.  Our result described above could be derived, with some work, from 
Theorem~5.2 combined with the information from the Remark~5.3 in
\cite{Kato70}.  See Theorem~3.2 of \cite{Yajima2011} for a Hilbert space version of Kato's Remark~5.3 in \cite{Kato70}. Our main result, Theorem~\ref{main1} below, does not follow from Kato's work but it reduces to a theorem of Kato if $X$ is a Hilbert space and $A(t)^{*}=-A(t)$, see Theorem 3 of \cite{Kato53} \footnote{This result of Kato seems to have been overlooked by many later authors.}. 
To explore the necessity of our assumptions, we give new counterexamples to well-posedness that are of a very simple and transparent type. 
Last but not least, the present paper shows that the essence of Kato's work in
the uniformly convex case can be summarized in a short and simple proof that 
requires nothing but basic functional analysis and a rudimentary knowledge of semigroup theory.

\section{Results and Examples}

Let $X$ be a complex Banach space and let $A(t):\dom\subset X\to X$ for
$t\in [0,T]$ be a family of closed linear operators with a time-independent dense domain $\dom\subset X$.
A two-parameter family of linear operators $U(t,s)\in \LL(X)$, will be
called an \emph{evolution system for $A(t)$ on} $D$ if the
following conditions are satisfied:
\begin{itemize}
\item[(i)] $U(t,s)\dom\subset \dom$ and the map $t\mapsto
  U(t,s)y$ on $[0,T]$ is a continuously differentiable solution of \eqref{IVP}
 for any $y\in \dom$ and $s\in [0,T]$. 
\item[(ii)] $U(s,s)=1$ and $U(t,r)U(r,s)=U(t,s)$ for all $s,r,t\in [0,T]$.
\item[(iii)] $(t,s)\mapsto U(t,s)$ strongly continuous on $[0,T]\times
  [0,T]$.
\end{itemize}
Any two-parameter family of linear operators $U(t,s)\in \LL(X)$
satisfying (ii) and (iii) is called an \emph{evolution system}.

Existence of an evolution system $U(t,s)$ with the properties analogous to 
(i)-(iii) on the triangle $0\leq s\leq t\leq T$ is equivalent to
well-posedness in the classical sense of $C^1$-solutions
\cite{EngelNagel}, Proposition~VI.9.3. Our assumptions on $A(t)$ in Theorem~\ref{main1} will allow us to
construct $U(t,s)$ on the entire square $[0,T]\times [0,T]$ and this is essential for our proof.

For the reader's convenience we recall that a Banach space $X$ is
called uniformly convex if, given $\eps>0$ there exists $\delta>0$ such that any pair
of normalized vectors $x,y\in X$ with $\|(x+y)/2\|>1-\delta$ satisfies
$\|x-y\|<\eps$. Every Hilbert space and every $L^p(\R^n)$ with
$1<p<\infty$ is uniformly convex. Uniformly convex Banach spaces are
reflexive (Milman) and uniform convexity implies that weak convergence
$x_n\weak x$ turns into strong convergence as soon as  $\|x_n\|\to \|x\|$.   

As a final preparation we recall from \cite{Pazy, EngelNagel} that the norm of every strongly
continuous semigroup $e^{At}$, $t\geq 0$, satifies a bound of the form
$\|e^{At}\|\leq Me^{\omega t}$. It follows that $\sigma(A)\subset \{\Rea z\leq \omega\}$. If $M=1$ then the semigroup
is called quasi-contractive. 

\begin{theorem}\label{main1}
Let $X$ be a uniformly convex Banach space and let $A(t):\dom\subset X\to
X$ for each $t\in [0,T]$ be the generator of a strongly continuous group with 
\begin{equation}\label{quasi}
   \|e^{A(t)s}\|\leq e^{\omega |s|},\qquad s\in \R,
\end{equation}
where $\omega$ and the domain $\dom$ are independent of $t$. Suppose that
$t\mapsto A(t)\in \LL(Y,X)$ is continuous and of bounded variation, where $Y$ is the space $D$ endowed with the graph norm of $A(0)$.
Then there exists a unique evolution system $U(t,s)$ for $A(t)$ on $\dom$. 
\end{theorem}

\begin{remark}\label{rem1}
The regularity assumption on $t\mapsto A(t)\in \LL(Y,X)$ is clearly satisfied if this map is Lipschitz, which, by the principle of uniform boundedness, is equivalent to 
$t\mapsto A(t)y$ being Lipschitz for all $y\in D$. 
\end{remark}

\begin{remark}
Theorem \ref{main1} is false if the assumption of uniform convexity is dropped (Example 1), and 
moreover, even if $X$ is a Hilbert space and $A(t)$ is skew-selfadjoint, the Lipschitz continuity cannot be replaced by H\"older continuity of some degree $\alpha<1$ (Example 2). This is in sharp contrast to the case of parabolic evolution equations, where H\"older continuity is sufficient \cite{Pazy, PruSch}.
\end{remark}

\begin{proof}
Let $\uA(t) :=A(t)-(\omega+1)$ and define $\|y\|_t :=\|\uA(t)y\|$ for $y\in D$. This norm is equivalent to the graph norm of $\uA(t)$ and hence $Y_t = (Y,\|\cdot\|_t)$ is 
a Banach space. Like $X$, the space $Y_t$ is uniformly convex as can be easily verified using the 
definition of uniform convexity given above. From now on $Y=Y_0$, which amounts to a different but equivalent choice of norm, compared to the definition of $Y$ in the theorem.

In the special case where $X$ is a Hilbert space and $A(t)$ is skew-selfadjoint, it follows that 
$\|y\|_t^2 = \|A(t)y\|^2+\|y\|^2$ and hence $Y_t$ is a Hilbert space too.

For $s,t\in [0,T]$ with $s<t$ we set 
\begin{equation}\label{def-V}
     V(s,t):=C\sup \sum_{i=1}^m \|A(t_i) - A(t_{i-1})\|_{Y,X}
\end{equation}
where the supremum is taken with respect to all partitions $s= t_0<t_1\ldots <t_m=t$ of the interval $[s,t]$. Up to the constant $C>0$, which will be chosen later, $V(s,t)$ is the variation of $\tau\mapsto A(\tau)$ over the interval $[s,t]$. Let $V(t,s):=V(s,t)$. Apart from the obvious inequality $V(t,s)\geq C\|A(t) - A(s)\|_{Y,X}$, the properties of $V$ that will be used in the following are, first, that 
\begin{equation}\label{add-V}
      V(r,s) + V(s,t) = V(s,t)\qquad \text{if}\ s<r<t,
\end{equation}
and, second, that $V(s,t)$ is separately continuous in $s$ and $t$. This follows from \eqref{add-V}, from the monotonicity of $V$, and from the continuity of $t\mapsto A(t)\in \LL(Y,X)$. The reader mainly interested in the case where $t\mapsto A(t)\in \LL(Y,X)$ is Lipschitz with some constant $L$ may replace $V(t,s)$ if $t>s$ by $CL(t-s)$ in all the following.

\medskip
\noindent\textbf{Step 1}: The constant $C$ in \eqref{def-V} may be chosen in such a way that for all
$s,t\in [0,T]$ and all $y\in\dom$,
$$
        \|y\|_t \leq e^{V(t,s)}\|y\|_s.
$$
\emph{Proof.} By the continuity of $t\mapsto \uA(t)\in\LL(Y,X)$, the map  $t\mapsto
\uA(t)^{-1}\in\LL(X,Y)$ is continuous and hence $C:=\sup_{s\in
  I}\|\uA(s)^{-1}\|_{X,Y}<\infty$. In view of $\|y\|_t\leq
\|\uA(t)\uA(s)^{-1}\| \|y\|_s$, Step 1 follows from 
\begin{align*}
    \|\uA(t)\uA(s)^{-1}\| &=  \|1+(\uA(t)-\uA(s))\uA(s)^{-1}\|\\
    &\leq 1+ \|\uA(t)-\uA(s)\|_{Y,X} \|\uA(s)^{-1}\|_{X,Y} \leq
    1+V(t,s) \leq e^{V(t,s)}. 
\end{align*}

\medskip
We now choose a sequence of partitions $\pi_n$ of $[0,T]$ with the
property that the mesh size  of $\pi_n$ vanishes in the limit
$n\to\infty$. Given $t\in [0,T]$ and $n\in \N$ we use $t_n$ to denote the largest element of $\pi_n$ less than or equal to $t$. 
The smallest element of $\pi_n$ larger than $t_n$ is denoted $t_n^{+}$, the largest one smaller than $t_n$ is denoted 
$t_n^{-}$. We thus have $t_n^{-}<t_n<t_n^{+}$ and 
$$
     t_n \leq t <t_n^{+}.
$$ 
Note that the points $t_n$ and $t_n^{\pm}$ are functions of both $t$ and $n$.
We define $U_n(t,s)$ for $t>s$ by 
$$
     U_n(t,s) := e^{A(t_n)(t-t_n)} e^{A(t_n^{-})(t_n-t_n^{-})}\cdots e^{A(s_n)(s_n^{+}-s)}  
$$
and $U_n(s,t):=U_n(t,s)^{-1}$. Note that $\|U_n(t,s)\|\leq
e^{\omega|t-s|}$ by assumption \eqref{quasi}.

\medskip
\noindent\textbf{Step 2}: For all $t>s$, $n\in\N$, and $y\in \dom$,
\begin{align*}
     \|U_n(t,s)y\|_t &\leq e^{V(t,s)+2V(s,s_n)+\omega(t-s)}\|y\|_s\\ 
     \text{and} \quad \|U_n(s,t)y\|_s &\leq e^{V(t,s)+2V(s,s_n)+\omega(t-s)}\|y\|_t.
\end{align*}
In particular, $\|U_n(t,s)\|_{Y,Y}<M$  for all $s,t\in [0,T]$ and all $n\in\N$. 

\smallskip\noindent
\emph{Proof.} 
With the help of Step 1 we pass from $\|\cdot \|_t$ to
$\|\cdot\|_{t_{n}}$, then from $\|\cdot\|_{t_{n}}$ to
$\|\cdot\|_{t_{n}^{-}}$ and so on, where in each step we use that
$e^{A(t)\tau}$ is a quasi-contraction in $Y_t$ satisfying
\eqref{quasi} for any $t\in [0,T]$. In this way we arrive at 
$$
        \|U_n(t,s)y\|_t \leq e^{V(t,s_n)+\omega(t-s)}\|y\|_{s_n},
$$
which, using Step 1 again, leads to the first of the asserted
inequalities. The second one is proved analogously and the uniform bound on $\|U_n(t,s)\|_{Y,Y}$ now follows from 
Step~1 and the compactness of $[0,T]$.

\medskip
\noindent\textbf{Step 3}: For all $x\in X$, the limit $U(t,s)x:=\lim_{n\to\infty} U_n(t,s)x$ exists
uniformly in $s,t\in [0,T]$. It defines an evolution system $U(t,s)$.

\smallskip\noindent
\emph{Proof.}  For any $y\in Y$ the map $\tau\mapsto U_m(t,\tau)
U_n(\tau,s)y$ is piecewise continuously differentiable with possible
jumps in the derivative at the partition points from $\pi_m\cup\pi_n$. It follows that
\begin{align*}
U_n(t,s)y - U_m(t,s)y  &= U_m(t,\tau) U_n(\tau,s)y |_{\tau=s}^{\tau=t}\\
 &=  \int_s^t U_m(t,\tau)\big(A(\tau_n)-A(\tau_m)\big) U_n(\tau,s)y\, d\tau.
\end{align*}
By Step 2 we conclude
$$
  \|U_n(t,s)y - U_m(t,s)y \| \leq \int_0^1 e^{\omega |t-\tau|}\|A(\tau_n) -
  A(\tau_m)\|_{Y,X} M\|y\|_Y \,d\tau \to 0\qquad (n,m\to\infty)
$$
by the continuity of $\tau\mapsto A(\tau)\in\LL(Y,X)$.
The assertion now follows from the density of $Y\subset X$ and from the uniform boundedness $\|U_n(t,s)\|\leq e^{\omega|t-s|}$. It follows that $(t,s)\mapsto U(t,s)x$ is continuous and the property (ii) of evolution systems is inherited from $U_n(t,s)$ as well.

\medskip
\noindent\textbf{Step 4}: $U(t,s)\dom\subset \dom$, and for all $y\in
D$ and $s,t\in [0,T]$,
$$
     \|U(t,s)y\|_t \leq e^{V(t,s)+\omega|t-s|}\|y\|_s.
$$
\noindent
\emph{Proof.} Let $y\in \dom$.  By Step~2 the sequence $(U_n(t,s)y)$
is bounded in $Y_t$ and by Step~3,  $U_n(t,s)y\to U(t,s)y$ in
$X$. Since $Y_t$ is reflexive it follows that $U(t,s)y\in D$ and that
$U_n(t,s)y\to U(t,s)y$ weakly in $Y_t$. Therefore, by the estimates of
Step~2,
$$
   \|U(t,s)y\|_t\leq \liminf_{n\to\infty} \|U_n(t,s)y\|_t\leq e^{V(t,s)+\omega |t-s|}\|y\|_s,
$$
where $V(s,s_n)\to 0$ as $n\to\infty$ was used.

\medskip
\noindent\textbf{Step 5}: For all $y\in \dom$ the map $t\mapsto U(t,s)y$
is differentiable in the norm of $X$ and 
$$
       \frac{d}{dt} U(t,s)y = A(t)U(t,s)y.
$$
\noindent
\emph{Proof.}  In view of $U(t,s)Y\subset Y$ and $U(t+h,s)=U(t+h,t)U(t,s)$, see Step 4, it suffices to prove the assertion for $s=t$. For any $y\in \dom$,
\begin{align*}
    U(t+h,t)y - e^{A(t)h} y &= \lim_{n\to\infty}
    e^{A(t)(h+t-\tau)} U_n(\tau,t)y \Big|_{\tau=t}^{\tau=t+h}\\
    &= \lim_{n\to\infty}\int_t^{t+h} e^{A(t)(h+t-\tau)} \big(A(\tau_n)-A(t)\big) U_n(\tau,s)y\, d\tau.
\end{align*}
By Step 2 it thus follows that 
\begin{eqnarray*}
\lefteqn{\Big\|\frac{1}{h} (U(t+h,t)y - e^{A(t)h} y)\Big\|}\\
 &\leq &\lim_{n\to\infty}\frac{1}{|h|}\left|\int_t^{t+h}
e^{\omega|t+h-\tau|} \|A(\tau_n)-A(t)\|_{Y,X}\, d\tau \right|\,
M\|y\|_Y\\
&= &\frac{1}{|h|}\left|\int_t^{t+h}
e^{\omega|t+h-\tau|} \|A(\tau)-A(t)\|_{Y,X}\, d\tau \right|\,
M\|y\|_Y \to 0\qquad (h\to 0)
\end{eqnarray*}
by the continuity of $\tau\mapsto A(\tau)\in \LL(Y,X)$. Since $(e^{A(t)h}y-y)/h\to A(t)y$ as $h\to 0$, 
Step 5 now follows.

\medskip
\noindent\textbf{Step 6}: For all $y\in D$ the map $t\mapsto A(t)U(t,s)y$ is continuous in the norm
of $X$.

\smallskip\noindent
\emph{Proof.} By the continuity of $t\mapsto A(t)\in \LL(Y,X)$ it
suffices to show that $t\mapsto U(t,s)y$ is continuous in the norm
of $Y$. To this end it suffices to show that $\lim_{h\to 0}U(t+h,t)y=
y$ in the norm of $Y$ or, equivalently, in the norm of $Y_t$. 
Since $U(t+h,t)y\to y$ in $X$ and since $h\mapsto U(t+h,t)y$ is
bounded in $Y_t$, see Step~1 and Step~4, it follows that $U(t+h,t)y\to
y$ weakly in $Y_t$. See the proof of Step~4 for a similar argument. Therefore,
\begin{align*}
    \|y\|_t &\leq \liminf_{h\to 0}\|U(t+h,t)y\|_t \leq \limsup_{h\to 0} \|U(t+h,t)y\|_t \\
             &\leq \limsup_{h\to 0}e^{V(t+h,t)}\|U(t+h,t)y\|_{t+h}
    \leq\limsup_{h\to 0} e^{2V(t+h,t)+\omega |h|}\|y\|_t = \|y\|_t. 
\end{align*}
The weak convergence $U(t+h,t)y\to y$ in $Y_t$ and the convergence of
the norms implies norm convergence in $Y_t$ by the uniform convexity.
\end{proof}

\begin{remark}\label{improved}
\begin{enumerate}
\item At the end of the introduction we pointed out two results of Kato that are closely related to Theorem \ref{main1}. There are two further prominent results in the literature on well-posedness, both due to Kato again, that can be compared to Theorem~\ref{main1}:
By a simple corollary of Theorem 1 of \cite{Kato73}, see Theorem 2.1.9 in \cite{Schmid-Diss}, it suffices to assume that 
\begin{equation}\label{A-map}
 t\mapsto A(t)\in\LL(Y,X)
\end{equation}
satisfies a certain $W^{1,1}_{*}$-regularity condition. This condition implies that \eqref{A-map} is absolutely continuous and hence continuous and of bounded variation. The condition $\partial_tA(t)\in L_{*}^{\infty}([0,T],\LL(Y,X))$ used by Kato in \cite{Kato85}, implies that \eqref{A-map} is Lipschitz
and hence continuous and of bounded variation.

\item If $A(t)$ was assumed to be the generator of a semigroup,
rather than a group, in Theorem~\ref{main1}, then the arguments of
our proof above still establish existence of a (unique) evolution system
$U(t,s)$ defined on the triangle $0\leq s\leq t\leq T$ such that 
$$
     \partial^{+}_tU(t,s)y = A(t)U(t,s)y,
$$
where $t\mapsto A(t)U(t,s)y$ is right-continuous and $\partial^{+}_t$
denotes the derivative from the right. Moreover, $ \partial_t U(t,s)y
= A(t)U(t,s)y$ except possibly for a countable set of $t$-values
depending on $y$ and $s$ (see the proof of Theorem~5.2 of \cite{Kato70}). 

\item In the case where the first or higher, suitably defined
  commutators of the operators $A(t)$ at distinct times are scalars, the continuity of the map
  \eqref{curve}, along with strong continuity of the commutators, is sufficient for well-posedness \cite{Goldstein,
    NickSch1998,Schmid2014}.  

\item  In the case where $X$ is a Hilbert space there are formal similarities between our Theorem~\ref{main1} 
  and the Theorem C.2 of Ammari and Breteaux \cite{AmmBre12}. In
  \cite{AmmBre12} the case of skew-selfadjoint generators with
  time-independent \emph{form domains} is considered and a notion
  of well-posedness in a weak sense is established.
\end{enumerate}
\end{remark}


In the remainder of this paper we specialize to operator families of
the form $A(t) = A_0 +B(t)$ where $A_0$ is the generator of a
$C_0$-group in $X$, $B(t)\in \LL(X)$ and $t\mapsto B(t)$ is strongly
continuous. Suppose that the evolution system $U(t,s)$ for $A$
exists. Then, for all $y\in D(A_0)$, 
$$
     U(t,s)y = e^{A_0(t-s)}y + \int_s^t d\tau e^{A_0(t-\tau)}B(\tau)U(\tau,s)y
$$
which, by assumption on $B(t)$, may be iterated indefinitely into a convergent Dyson
series \cite{ReedSimon}. For the evolution system in the interaction picture we obtain
\begin{align}
   e^{-A_0 t}U(t,s)e^{A_0 s}y = y+ &\int_s^t d\tau_1 \tilde B(\tau_1)y\nonumber\\
               +&\int_s^t  d\tau_1 \int_s^{\tau_1}d\tau_2 \tilde
               B(\tau_1)\tilde B(\tau_2)y + \ldots,\label{Dyson}
\end{align}
where $\tilde B(\tau) y := e^{-A_0\tau} B(\tau) e^{A_0\tau}$. Consequently, the operator family $U(t,s)$ defined by \eqref{Dyson} is
the only candidate for the evolution system generated by
$A(t)=A_0+B(t)$. 

The following theorem is now an immediate  corollary of the previous
one and Theorem~3.1.1 from \cite{Pazy}. In the uniformly convex case it improves on a similar result due to Phillips: in Theorem~6.2 of \cite{Phillips1953} it is assumed that $t\mapsto B(t)$ is strongly continuously differentiable.

\begin{theorem}\label{main2}
Suppose $X$ is a uniformly convex Banach space and that $A(t)=A_0+B(t)$ for $t\in
[0,T]$, where $A_0:D\subset X\to X$ is the generator of a strongly
continuous quasi-contractive group in $X$ and $B(t)\in\LL(X)$. If $t\mapsto B(t)y$ is
continuous for all $y\in X$ and Lipschitz for all $y\in D$, then there exists a unique evolution
system $U(t,s)$ for $A$ on $D$ and $e^{-A_0 t}U(t,s)e^{A_0 s}$ is given by the Dyson series \eqref{Dyson}.
\end{theorem}


The following examples show that the assumptions of  uniform convexity
and Lipschitz continuity in this theorem and hence in the
Theorem~\ref{main1} cannot be weakened in an essential way.

\medskip\noindent
\textbf{Example 1.} Let $X=C_0(\R)$ be the Banach space of bounded and
continuous functions vanishing at infinity, the norm being the usual
maximum norm. This space is not uniformly convex. Let  $e^{A_0t}$ be the
strongly continuous group in $X$ defined by left translations, that is, $e^{A_0t}x(\xi)=x(\xi+t)$. We define 
$A(t):D\subset X\to X$ for $t\in [0,1]$ by $D=D(A_0)$ and 
$$
     A(t)= A_0 + B(t),\qquad B(t)=e^{A_0t}Be^{-A_0t}
$$
where $B$ denotes multiplication with the following bounded function $f:\R\to [0,1]$: we choose $f(\xi)=0$ for $\xi\leq 0$, $f(\xi)=\xi$ for $\xi\in [0,1]$ and
$f(\xi)=1$ for $\xi\geq 1$.  Then $B(t)$ is multiplication with the
function $\xi\mapsto f(\xi+t)$ and from the fact that $f$ is Lipschitz it is
easy to check that $t\mapsto B(t)$ is strongly Lipschitz. If an
evolution system $U(t,s)$ for $A$ on $D$ existed, then it would be given by the Dyson
series \eqref{Dyson} and since $\tilde B(t) = B$ it would follow that 
\begin{equation}\label{dyson-evolution}
       U(t,0) = e^{A_0 t} e^{Bt}.
\end{equation}
Since $D=D(A_0)$ is left invariant by $e^{-A_0t}$ it would follow that $e^{Bt}D(A_0)\subset D(A_0)$.
But $D(A_0) = \{y\in C^1(\R) \mid y,y'\in X\}$ and the operator $e^{Bt}$ acts as multiplication with the
non-differentiable function $e^{f(\xi)t}$. Hence
$e^{Bt}D(A_0)\not\subset D(A_0)$ and we have a
contradiction. Therefore an evolutions system $U$ for $A$ on $D$ cannot exist.

\medskip\noindent
\textbf{Example 2.} For this example we adopt all elements of Example~1 with two exceptions: now $X=L^2(\R)$ and 
$f$ denotes multiplication with the bounded function $f=ig$, where $g:\R\to\R$ is the Weierstra{\ss} function
$$
     g(\xi) = \sum_{n=1}^{\infty} 2^{-n}\cos(2^n\xi).
$$
This function is H\"older continuous of degree $\alpha$ for all $\alpha<1$ and nowhere differentiable. See \cite{Zygmund}, Theorem~II.4.9, including the proof, and the remark after Theorem~II.4.10. It easily follows that $t\mapsto B(t)$ is strongly H\"older continuous of degree $\alpha$
for all $\alpha<1$. As in Example~1 we argue that $e^{Bt}D(A_0)\subset D(A_0)$ if the evolution system $U$ for $A$ existed.
But $e^{Bt}$ acts by multiplication with $\xi\mapsto e^{f(\xi)t}$, which is
nowhere differentiable, and $D(A_0) = H^1(\R)$ whose elements are differentiable almost everywhere. We have a contradiction 
and hence an evolution system $U$ for $A$ on $D$ cannot exist. Note that $A(t)$ is skew-selfadjoint in this example.

\medskip\noindent
\textbf{Acknowledgement:} J. S. would like to thank Roland Schnaubelt for interesting discussions and for pointing out \cite{PruSch} to us. His work was supported by the \emph{Deutsche Forschungsgemeinschaft (DFG)} through the Research Training Group 1838.

\vspace{\baselineskip}

\end{document}